\documentclass[12pt,oneside,reqno]{amsart}
\usepackage{txfonts}
\usepackage{bbm}
\usepackage{amsmath}
\usepackage{graphicx}
\usepackage{mathrsfs}
\usepackage{stmaryrd}
\usepackage{amsfonts}
\usepackage{enumerate,amsmath,amssymb,amsthm}

\numberwithin{equation}{section}

\newcommand{\be}{\begin{eqnarray}}
\newcommand{\ee}{\end{eqnarray}}
\newcommand{\ce}{\begin{eqnarray*}}
\newcommand{\de}{\end{eqnarray*}}
\newtheorem{theorem}{Theorem}[section]
\newtheorem{lemma}[theorem]{Lemma}
\newtheorem{remark}[theorem]{Remark}
\newtheorem{definition}[theorem]{Definition}
\newtheorem{proposition}[theorem]{Proposition}
\newtheorem{Examples}[theorem]{Example}
\newtheorem{corollary}[theorem]{Corollary}

\def\eps{\varepsilon}

\def\b{\beta}
\def\p{\partial}

\def\[{{\Big[}}
\def\]{{\Big]}}
\def\<{{\langle}}
\def\>{{\rangle}}
\def\({{\Big(}}
\def\){{\Big)}}

\def\bx{{\mathbf{x}}}
\def\tr{{\rm tr}}

\def\dif{{\mathord{{\rm d}}}}

\def\no{\nonumber}
\def\={&\!\!=\!\!&}
\def\bt{\begin{theorem}}
\def\et{\end{theorem}}
\def\bl{\begin{lemma}}
\def\el{\end{lemma}}
\def\br{\begin{remark}}
\def\er{\end{remark}}

\def\bd{\begin{definition}}
\def\ed{\end{definition}}
\def\bp{\begin{proposition}}
\def\ep{\end{proposition}}
\def\bc{\begin{corollary}}
\def\ec{\end{corollary}}
\def\bx{\begin{Examples}}
\def\ex{\end{Examples}}

\def\cB{{\mathcal B}}

\def\cF{{\mathcal F}}

\def\cL{{\mathcal L}}
\def\cM{{\mathcal M}}

\def\cP{{\mathcal P}}

\def\cT{{\mathcal T}}

\def\mD{{\mathbb D}}
\def\mE{{\mathbb E}}

\def\mI{{\mathbb I}}

\def\mN{{\mathbb N}}

\def\mR{{\mathbb R}}

\def\mT{{\mathbb T}}

\def\mW{{\mathbb W}}

\def\bP{{\mathbf P}}

\def\sF{{\mathscr F}}

\def\geq{\geqslant}
\def\leq{\leqslant}

\def\div{\mathord{{\rm div}}}

\def\bP{{\mathbf P}}

\def\tr{{\mathrm t}}
\allowdisplaybreaks

\begin{document}

\title{Stochastic Lagrangian Particle Approach to Fractal Navier-Stokes Equations}

\date{}
\author{Xicheng Zhang}

%\thanks{$*$ This work is supported by ARC Discovery grant DP0663153 of Australia.}
%Galerkin's Approximation, Analytic Semigroup,

\dedicatory{
School of Mathematics and Statistics,
Wuhan University, Wuhan, Hubei 430072, P.R.China,\\
Email: XichengZhang@gmail.com
 }

\begin{abstract}
In this article we study the fractal Navier-Stokes equations by
using stochastic Lagrangian particle path approach in Constantin and Iyer \cite{Co-Iy}.
More precisely, a stochastic representation for the fractal Navier-Stokes equations is given in terms of
stochastic differential equations driven by L\'evy processes. Basing on this representation, a self-contained proof for the
existence of local unique solution for the fractal Navier-Stokes equation with initial data in $\mW^{1,p}$ is provided,
and in the case of two dimensions or large viscosity, the existence of global solution is also obtained. In order to obtain the global existence in any dimensions for large viscosity,
the gradient estimates for L\'evy processes with time dependent and discontinuous drifts is proved.
\end{abstract}

\maketitle
\rm

\section{Introduction}
Consider the following incompressible fractal or generalized Navier-Stokes equation in $\mR^d$ (abbreviated as FNSE):
\be
\left\{
\begin{aligned}
&\p_t u=\cL u-(u\cdot\nabla)u+\nabla p,\ \ t\geq 0,\\
&\nabla\cdot u=0,\ \ u(0)=u_0,
\end{aligned}
\right.\label{FNS0}
\ee
where $u=(u^1,\cdots, u^d)^{\mathrm{t}}$ denotes the column vector of velocity field, $p$ is the pressure, $\cL$ is the generator of a L\'evy process given by
\begin{align}
\cL u(x)=\int_{\mR^d} (u(x+y)-u(x)-1_{|y|\leq 1}(y\cdot\nabla)u(x))\nu(\dif y),\label{EP}
\end{align}
where $\nu$ is a L\'evy  measure on $\mR^d$, i.e., it satisfies that $\nu\{0\}=0$ and
$$
\int_{\mR^d}1\wedge|y|^2\nu(\dif y)<+\infty.
$$
When $\nu(\dif y)=\dif y/|y|^{d+\alpha}$ with $\alpha\in(0,2)$, $\cL=-c_\alpha(-\Delta)^{\alpha/2}$ is the usual fractional Laplacian operator by
multiplying a constant.

As a simplified model of equation (\ref{FNS0}), the following fractal Burgers equation has been studied by
Biler, Funaki and Woyczynski \cite{Bi-Fu-Wo} and Kiselev, Nazarov and Schterenberg \cite{Ki-Na-Sh},
$$
\p_t u=-(-\Delta)^{\alpha/2} u-(u\cdot\nabla)u,\ \ t\geq 0,\ \ u(0)=u_0.
$$
As for generalized Navier-Stokes equation (\ref{FNS0}), when $\cL=-(-\Delta)^{\alpha/2}$, it has been studied by Wu \cite{Wu}  in Besov spaces by using purely analytic argument.
The main feature of such fractal equations is that operator $\cL$ given by (\ref{EP}) is non-local. Recently, there are increasing interests
for studying such fractal equations or fractional dissipative equations since they naturally appear
in hydrodynamics, statistcal mechanics, physiology, certain combustion models, and so on (cf. \cite{Sh-Za-Fr, Re, Wo}, etc.).

The aim of this paper is to study equation (\ref{FNS0}) by using a stochastic Lagrangian particle trajectories approach
following \cite{Co-Iy, Zh2}. More precisely,  Constantin and Iyer \cite{Co-Iy} gave the following elegant stochastic representation
for the regularity solution $u$ of Navier-Stokes equation (corresponding to $\cL=\nu\Delta$ in (\ref{FNS0})):
\begin{align}
\label{Re}\left\{
\begin{aligned}
&X_t(x)=x+\int^t_0u_s(X_s(x))\dif s+\sqrt{2\nu}B_t,\ \ t\geq 0,\\
&u_t=\bP\mE[(\nabla^\tr X^{-1}_t) (u_0\circ X^{-1}_t)],
\end{aligned}
\right.
\end{align}
where $\bP$ denotes the Leray-Hodge projection onto divergence free vector fields, $B_t$ is a Brownian motion,
and $X^{-1}_t(x)$ denotes the inverse of $x\mapsto X_t(x)$.
Basing on this representation, a self-contained proof of the existence of local smooth solutions in H\"older space was given by Iyer \cite{Iy0}.
Later on, by reversing the time variable, in a previous work \cite{Zh2}, we considered the following stochastic representation:
\begin{align}
\label{SS}\left\{
\begin{aligned}
&X_{t,s}(x)=x+\int^s_tu_r(X_{t,r}(x))\dif r+\sqrt{2\nu}(B_s-B_t),\ \ t\leq s\leq 0,\\
&u_t=\bP\mE[(\nabla^\tr X_{t,0}) (u_0\circ X_{t,0})],
\end{aligned}
\right.
\end{align}
and a self-contained proof of the existence of local smooth solutions in Sobolev space is also obtained.
Moreover, the global solution for large viscosity is proven by using Bismut formula.

Naturally, if one replaces the Brownian motion in (\ref{SS}) by a general L\'evy process $L_t$, then it is expected that the corresponding
solution $u$ will solve the following backward fractal Navier-Stokes equation:
\begin{align}
\left\{
\begin{aligned}
&\p_t u+\cL u+(u\cdot\nabla)u+\nabla p=0,\ \ t\leq 0,\\
&\nabla\cdot u=0,\ \ u(0)=u_0,
\end{aligned}
\right.\label{FNS}
\end{align}
where $\cL$ is the L\'evy generator of the L\'evy process $L_t$. In general, it seems hard to solve the above fractal Navier-Stokes equation
by using purely analytic tools. However, stochastic system (\ref{SS}) is easier to be dealt with if one replaces $B_t$
by a general process and only considers the local smooth solutions. In fact, it is easy to obtain the existence of local smooth
solutions for stochastic system (\ref{SS}), and a global smooth solution in two dimensional case by the same arguments as in \cite{Iy0, Zh3}.
This will be given in Section 2.

On the other hand, if we only assume that the initial data belongs to the first order Sobolev space $\mW^{1,p}$, it seems not so easy to construct
a local solution in $\mW^{1,p}$ for stochastic system (\ref{SS}).  A clear difficulty is to obtain the differentiability of the solution flow $x\mapsto X_{t,0}(x)$.
Although one can solve the following equation with  divergence free vector field $u\in L^1_{loc}((-\infty,0];\mW^{1,p})$ by using DiPerna-Lions' theory (cf. \cite{Di-Li, Cr-De-Le}):
$$
X_{t,s}(x)=x+\int^s_tu_r(X_{t,r}(x))\dif r+(L_s-L_t),
$$
it is only known that $x\mapsto X_{t,0}(x)$ is approximately differentiable (cf. \cite{Am-Le-Ma, Cr-De-Le}). This difficulty will be overcomed by using Krylov's estimate for jump-diffusion
processes and the regularizing effect of L\'evy process if $L_t$ is non-degenerated in some sense (see Condition {\bf (H)$_\alpha$} below).
Thus, following Section 2, Section 3 will be devoted to the proof of the existence of local $\mW^{1,p}$-solutions.

For proving the existence of global solutions in $\mW^{1,p}$ for large viscosity, we need some gradient estimates for the SDE with Sobolev coefficients and driven by a L\'evy process.
For this aim, we shall use some asymptotic estimates for the heat kernels of L\'evy processes due to Schilling, Sztonyk and  Wang \cite{Sc-Sz-Wa}. Our approach for the gradient estimates of SDEs
is based on the a priori estimate for an integro-differential equation and the uniqueness of weak solutions. This is the content of Section 4, and can be read independently.
In Section 5, we prove the global well posedness for large viscosity.

Lastly, we mention that other stochastic approaches for incompressible Navier-Stokes equations can be found in the references \cite{Le-Sz, Bu, Bu-Fl-Ro,Cr-Sh,Es-Ma-Pu-Sc}, etc.;
and compared with the analytic arguments, one of the main advantages of representation (\ref{SS}) is that it is convenient for numerical simulations (cf. \cite{Ma-Be,Iy-Ma}).
This is in fact our main motivation for studying the stochastic representation of fractal Navier-Stokes equation (\ref{FNS}).

\section{Stochastic representation for fractal Navier-Stokes equations}

We first fix some notations. Set $\mN_0:=\{0\}\cup\mN$ and $\mR_-:=(-\infty,0]$. For $k\in\mN_0$, let
$C^k_b=C^k_b(\mR^d;\mR^d)$ be the space of all $k$-order continuously differentiable vector fields on $\mR^d$ with the norm
$$
\|u\|_{C^k_b}:=\sum_{j=0}^k\sup_{x\in\mR^d}|\nabla^j u(x)|<+\infty,
$$
where $\nabla^j$ denotes the $j$-order gradient, and $|\cdot|$ denotes the Euclidean norm.
For $k\in\mN_0$ and $p\geq 1$, let $\mW^{k,p}=\mW^{k,p}(\mR^d;\mR^d)$ be the usual vector-valued Sobolev space on $\mR^d$ with the norm
$$
\|u\|_{k,p}:=\sum_{j=0}^k\|\nabla^ju\|_p<+\infty,
$$
where $\|\cdot\|_p$ is the usual $L^p$-norm in $\mR^d$.

Let us now recall some basic notions and facts about L\'evy processes on negative time axis.
Let $(L_t)_{t\leq 0}$ be an $\mR^d$-valued L\'evy  process on some probability space $(\Omega,\sF,P)$, i.e.,
an $\mR^d$-valued stochastically continuous  process with stationary independent increments and $L_0=0$.
By L\'evy-Khintchine's formula (cf. \cite[p.109, Corollary 2.4.20]{Ap}), the characteristic function of $L_t$ is given by
\begin{align}
\mE(e^{i\xi\cdot L_t})
=\exp\left\{t\left[ib\cdot\xi+\xi^{\mathrm{t}} A\xi+\int_{\mR^d}[1-e^{i\xi\cdot x}+i\xi\cdot x 1_{|x|\leq 1}]\nu(\dif x)\right]\right\}
=:e^{t\psi(\xi)},\label{Lp9}
\end{align}
where $\psi(\xi)$ is a complex-valued function called the symbol of $(L_t)_{t\leq 0}$, and
$b\in\mR^d$, $A\in\mR^d\times\mR^d$ is a positive definite and symmetric matrix,
$\nu$ is a L\'evy  measure on $\mR^d$. Throughout this paper, we only consider the pure jump L\'evy process, and assume below that
$$
b=0,\ \ A=0.
$$
We remark that $t\mapsto L_t$ admits a version still denoted by $L_t$ such that for almost all $\omega$,
$t\mapsto L_t(\omega)$ is right continuous and has left limit, but,
for fixed $t$,
$$
P\{\omega: L_t(\omega)\not= L_{t-}(\omega)\}=0.
$$
Below, for $t\leq s\leq 0$, define
\begin{align}
\sF_{t,s}:=\sigma\{L_r-L_t: t\leq r\leq s\}.\label{Flow}
\end{align}

Given $u\in C(\mR_{-}; C^3_b(\mR^d;\mR^d))$, for $x\in\mR^d$, let $X_{t,s}(x)$ be the unique solution of
the following SDE:
\begin{align}
X_{t,s}(x)=x+\int^s_t u_r(X_{t,r}(x))\dif r+(L_s-L_t),\ \ t\leq s\leq 0.\label{Eq1}
\end{align}
It is easy to see that $\{X_{t,s}(x), x\in\mR^d, t\leq s\leq 0\}$ forms a stochastic $C^3$-diffeomorfism flow, and
\begin{align}
\nabla X_{t,s}(x)=\mI+\int^s_t\nabla u_r(X_{t,r}(x))\cdot\nabla X_{t,r}(x)\dif r,\label{Lp10}
\end{align}
where $\nabla X_{t,s}(x)=(\p_jX^i_{t,s}(x))_{i,j=1,\cdots, d}$, and $\p_j$ denotes the partial derivative with respect to the $j$-th component of $x$.

Let $N(t,\Gamma):=\sum_{t\leq s<0}1_{\Gamma}(L_s-L_{s-}), \Gamma\in\cB(\mR^d)$ be the Poisson random point measure associated
with $(L_t)_{t\leq 0}$. By L\'evy-It\^o's decomposition (cf. \cite[p.108, Theorem 2.4.16]{Ap}),  one has
$$
L_t=\int_{|x|\leq 1}x\tilde N(t,\dif x)+\int_{|x|>1}xN(t,\dif x),
$$
where $\tilde N(t,\dif x):=N(t,\dif x)-t\nu(\dif x)$ is the compensated random martingale measure.
For $g\in C^2_b(\mR^d;\mR)$, by It\^o's formula (cf. \cite[p.226, Theorem 4.4.7]{Ap}), we have
\begin{align}
g(X_{t,s})=g(x)+\int^s_t [\cL g(X_{t,r})+(u_r\cdot\nabla) g(X_{t,r})]\dif r+M^g_{t,s},\label{Ito}
\end{align}
where $\cL$ is the generator of $(L_t)_{t\leq 0}$ given by (\ref{EP}),
and $(M^g_{t,s})_{s\in[t,0]}$ is a square integrable martingale given by
$$
M^g_{t,s}:=\int^s_t\!\!\!\int_{\mR^d}[g(X_{t,r-}+y)-g(X_{t,r-})]\tilde N(\dif r,\dif y).
$$

We have
\bt\label{Th2}
Let $u\in C(\mR_-;C^3_b(\mR^d;\mR^d))$ and $X_{t,s}(x)$ be the solution of SDE (\ref{Eq1}).
For $\varphi\in C^2_b(\mR^d;\mR^d)$ and $c\in C(\mR_-; C^2_b(\mR^d;\mR))$, define
$$
h_t(x):=\mE\left[\exp\left\{\int^0_tc_r(X_{t,r}(x))\dif r\right\}\varphi(X_{t,0}(x))\right],
$$
and
$$
w_t(x):=\mE[\nabla^{\mathrm{t}} X_{t,0}(x)\cdot\varphi(X_{t,0}(x))].
$$
Then $h, w\in C^1(\mR_-; C^2_b(\mR^d;\mR^d))$ respectively and uniquely solve the following partial integro-differential equations (PIDE):
\begin{align}
\p_t h_t+\cL h_t+(u_t\cdot\nabla)h_t+c_t h_t=0,\ \ h_0(x)=\varphi(x),\label{Eq2}
\end{align}
and
\begin{align}
\p_t w_t+\cL w_t+(u_t\cdot\nabla)w_t+(\nabla^{\mathrm{t}} u_t)w_t=0,\ \ w_0(x)=\varphi(x).\label{Eq3}
\end{align}
\et
\begin{proof}
Fix $t<0$. For $g\in C^2_b(\mR^d)$ and $\delta>0$, by It\^o's formula (see (\ref{Ito})), we have
\begin{align*}
&\mE\left[\exp\left\{\int^t_{t-\delta}c_r(X_{t-\delta,r}(x))\dif r\right\}g(X_{t-\delta,t}(x)) \right]-g(x)\\
&\quad=\mE\left[\int^{\mathrm{t}}_{t-\delta}\exp\left\{\int^s_{t-\delta}c_r(X_{t-\delta,r}(x))\dif r\right\}
\Big[\cL g(X_{t-\delta,s}(x))+(u_s\cdot\nabla)g(X_{t-\delta,s}(x))\Big]\dif s\right]\\
&\qquad+\mE\left[\int^t_{t-\delta}\exp\left\{\int^s_{t-\delta}c_r(X_{t-\delta,r}(x))\dif r\right\} c_s(X_{t-\delta,s}(x))g(X_{t-\delta,s}(x))\right]
\end{align*}
By the stochastic continuity of $t\mapsto L_t$,  from equation (\ref{Eq1}), it is easy to prove
that $(t,s)\to X_{t,s}(x)$ is also stochastically continuous. Thus,
since $(s,x)\mapsto \cL g(x)+(u_s\cdot\nabla)g(x)+c_s(x) g(x)$ is bounded and continuous, by the dominated convergence theorem, we have
\begin{align}
&\frac{1}{\delta}\left[\mE\left[ \exp\left\{\int^t_{t-\delta}c_r(X_{t-\delta,r}(x))\dif r\right\}g(X_{t-\delta,t}(x))\right]-g(x)\right]\no\\
&\quad\stackrel{\delta\to 0}{\to} \cL g(x)+(u_t\cdot\nabla)g(x)+c_t(x) g(x).\label{Lp01}
\end{align}
Noticing that
$$
X_{t-\delta,0}(x)=X_{t,0}\circ X_{t-\delta,t}(x),
$$
by Markov property, we have
\begin{align*}
h_{t-\delta}&=\mE\left[\exp\left\{\int^t_{t-\delta}c_r(X_{t-\delta,r})\dif r\right\}
\exp\left\{\int^0_tc_r(X_{t,r}\circ X_{t-\delta,t})\dif r\right\}\varphi(X_{t,0}\circ X_{t-\delta,t})\right]\\
&=\mE\left[\exp\left\{\int^t_{t-\delta}c_r(X_{t-\delta,r})\dif r\right\}
\mE\left[\exp\left\{\int^0_tc_r(X_{t,r}\circ X_{t-\delta,t})\dif r\right\}\varphi(X_{t,0}\circ X_{t-\delta,t})\Big|\sF_{t-\delta,t}\right]\right]\\
&=\mE\left[\exp\left\{\int^t_{t-\delta}c_r(X_{t-\delta,r})\dif r\right\}h_t\circ X_{t-\delta,t}\right].
\end{align*}
Thus, by (\ref{Lp01}), we obtain
\begin{align*}
\frac{1}{\delta}(h_t(x)-h_{t-\delta}(x))\stackrel{\delta\to 0}{\to} -[\cL h_t(x)+(u_t\cdot\nabla)h_t(x)+c_t(x) h_t(x)].
\end{align*}
Since the limit is a continuous function of $(t,x)$, it follows that for each $x$, $t\mapsto h_t(x)$ is differentiable and equation (\ref{Eq2})
is obtained.

As for (\ref{Eq3}), observing that
$$
\nabla X_{t-\delta,0}(x)=(\nabla X_{t,0})\circ X_{t-\delta,t}(x)\cdot\nabla X_{t-\delta,t}(x),
$$
by Markov property again, we have
\begin{align*}
w_{t-\delta}(x)=\mE[\nabla^{\mathrm{t}} X_{t-\delta,0}(x)\cdot\varphi(X_{t-\delta,0}(x))]
=\mE[\nabla^{\mathrm{t}} X_{t-\delta,t}(x)\cdot w_t(X_{t-\delta,t}(x))].
\end{align*}
Hence, by (\ref{Lp01}) and (\ref{Lp10}), we have
\begin{align*}
\frac{1}{\delta}(w_t(x)-w_{t-\delta}(x))&=-\frac{1}{\delta}\mE[w_t(X_{t-\delta,t}(x))-w_t(x)]
-\frac{1}{\delta}\mE[(\nabla^{\mathrm{t}} X_{t-\delta,t}-\mI)\cdot w_t(X_{t-\delta,t}(x))]\\
&\stackrel{\delta\to 0}{\to}-[\cL w_t+(u_t\cdot\nabla)w_t](x)-(\nabla^{\mathrm{t}} u_t\cdot w_t)(x).
\end{align*}
Equation (\ref{Eq3}) is thus obtained.

We now prove the uniqueness. Here we adopt the duality argument. Let $\hat X_{t,s}(x)$ solve the following SDE:
$$
\hat X_{t,s}(x)=x-\int^s_t u_r(\hat X_{t,r}(x))\dif r-(L_s-L_t),\ \ t\leq s\leq 0.
$$
Fix $\phi\in C^\infty_0(\mR^d;\mR^d)$ and $T<0$. For $t\in[T,0]$, define
$$
\hat h_t(x):=\mE\left[\exp\left\{\int^0_{T-t}(c_r-\div u_r)(X_{T-t,r}(x))\dif r\right\}\phi(X_{T-t}(x))\right].
$$
By the above proof, it follows that $\hat h_t\in L^1(\mR^d)\cap C^2_b(\mR^d)$ solves the following PIDE:
\begin{align*}
\p_t \hat h_t&=\cL^*\hat h_t-(u_t\cdot\nabla)\hat h_t+(c_t-\div u_t)\hat h_t\\
&=\cL^* \hat h_t-\div (u_t\otimes \hat h_t)+c_t \hat h_t
\end{align*}
subject to $\hat h_T(x)=\phi(x)$, where $\cL^*$ is the dual operator of $\cL$ and given by
$$
\cL^* g(x)=\int_{\mR^d} [g(x-y)-g(x)+(y\cdot \nabla )g(x)1_{|y|\leq 1}]\nu(\dif y).
$$
Now, let $h\in C^1(\mR_-; C^2_b(\mR^d;\mR^d))$ solve (\ref{Eq2}) with $h_0(x)\equiv 0$. Then,
by the integration by parts formula, we have
$$
\p_t\<h_t,\hat h_t\>=-\<\cL h_t+(u_t\cdot\nabla)h_t+c_t h_t,\hat h_t\>
+\<h_t,\cL^* \hat h_t-\div (u_t \otimes\hat h_t)+c_t \hat h_t\>=0,
$$
where $\<h_t,\hat h_t\>=\int_{\mR^d}\<h_t(x),\hat h_t(x)\>_{\mR^d}\dif x$.
Since $\<h_0,\hat h_0\>=0$, it is immediate that $\<h_T,\hat h_T\>=\<h_T,\phi\>=0$, which then gives $h_T(x)=0$ by the arbitrariness
of $\phi\in C^\infty_0(\mR^d;\mR^d)$.
\end{proof}

\br
If one assumes that $u\in L^1_{loc}(\mR_-;C^3_b(\mR^d;\mR^d))$ and $c\in L^1_{loc}(\mR_-; C^2_b(\mR^d;\mR))$, then the conclusions of Theorem \ref{Th2}
still hold if one replaces equations (\ref{Eq2}) and (\ref{Eq3}) by
$$
h_t(x)=\varphi(x)+\int^0_t[\cL h_s(x)+(u_s\cdot\nabla)h_s(x)+c_s(x) h_s(x)]\dif s,
$$
and
$$
w_t(x)=\varphi(x)+\int^0_t[\cL w_s(x)+(u_s\cdot\nabla)w_s(x)+(\nabla^{\mathrm{t}} u_s)w_s(x)\dif s.
$$
\er

Using Theorem \ref{Th2}, we have the following representation for the solution of fractal Navier-Stokes equation
as in \cite[Theorem 2.3]{Zh2}.
\bt\label{Rep}
Let $u\in C(\mR_-;C^3_b(\mR^d;\mR^d))$ be divergence free. Then, $u$ is a solution of fractal Navier-Stokes equation (\ref{FNS}) if and only if $u$
solves the following stochastic system:
\begin{align}\label{SN}
\left\{
\begin{aligned}
&X_{t,s}(x)=x+\int^s_t u_r(X_{t,r}(x))\dif r+(L_s-L_t),\ \ t\leq s\leq 0,\\
&u_t=\bP\mE[\nabla^{\mathrm{t}} X_{t,0}\cdot(u_0\circ X_{t,0})],
\end{aligned}
\right.
\end{align}
where $L$ is a L\'evy  process with generator $\cL$, and $\bP$ is the Leray-Hodge projection onto divergence free vector fields.
\et

Along the completely same lines as in \cite[Theorems 3.8 and 4.2]{Zh2}, one has the following result. The details are omitted.
\bt\label{Main}
For any $k\in\mN_0$ and $p>d$, there exists a constant $C_0=C_0(k,p,d)>0$ such that for
any $u_0\in \mW^{k+2,p}(\mR^d;\mR^d)$ divergence free and $T:=-(C_0\|\nabla u_0\|_{k+1,p})^{-1}$,
there is a unique pair of $(u,X)$ with $u\in C([T,0];\mW^{k+2,p})$ satisfying (\ref{SN}).
Moreover, for any $t\in[T,0]$,
\begin{align*}
\|\nabla u_t\|_{k+1,p}\leq C_0\|\nabla u_0\|_{k+1,p}.
\end{align*}
In two dimensional case, one has that for all $t\in\mR_-$,
$$
\|u_t\|_{k+2,p}\leq C(\|u_0\|_{k+2,p},k,p,t),
$$
where the constant $C$ is increasing with respect to its first argument.
In particular, there exists a unique global solution $u\in C(\mR_-;\mW^{k+2,p})$ to (\ref{SN}) in the two dimensional case.
\et

\section{Existence of local solutions for FNSE with $\mW^{1,p}$ initial data}

In the remaining sections, we mainly study equation (\ref{SN}) with $u_0\in\mW^{1,p}(\mR^d;\mR^d)$. For this aim, we  assume that
\begin{enumerate}[{\bf (H)$_\alpha$}]
\item Let $\psi(\xi)$ be the L\'evy symbol given in (\ref{Lp9}) and satisfy that for some $\alpha\in(0,2)$,
$$
\mathrm{Re}(\psi(\xi))\asymp |\xi|^\alpha\ \mbox{ as }\ |\xi|\to\infty,
$$
where $a\asymp b$ means that for some $c_1,c_2>0$, $c_1 b\leq a\leq c_2 b$.
\end{enumerate}

Consider the following SDE:
\begin{align}
X_{t,s}(x)=x+\int^s_t u_r(X_{t,r}(x))\dif r+\nu^{1/\alpha}(L_s-L_t),\ \ t\leq s\leq 0,\label{EL7}
\end{align}
where $u:\mR_-\times\mR^d\to\mR^d$ is a bounded Borel measurable function, and with a little abuse of notations,
$\nu\geq 0$ denotes a positive constant which plays the viscosity role.

We recall the following Krylov estimate for jump diffusion processes taken from \cite[Theorem 3.7]{Zh3}.
Although the theorem is given therein only for $\alpha$-stable processes, it is clearly also valid for more general L\'evy processes
considered in the present paper since the proof only depends on the gradient estimate (\ref{Pl1}) below.
\bt\label{Kry1}
Suppose  that {\bf (H)$_\alpha$} holds with $\alpha\in(1,2)$, and $u$ is bounded by $\kappa$. Let $X_{t,s}(x)$ solve equation (\ref{EL7}).
Then for any $p>\frac{d}{\alpha}$ and $q>\frac{p\alpha}{p\alpha-d}$,
there exists a constant $C_{\kappa}=C_{\kappa}(d,\alpha,p,q,\psi)>0$ independent of $\nu\geq 1$ and $x\in\mR^d$
such that for  all $-1\leq t\leq s_1<s_2\leq 0$  and  $f\in L^q([s_1,s_2]; L^p(\mR^d))$,
\begin{align}
\mE\left(\int^{s_2}_{s_1} f_r(X_r(x))\dif r\Big|\sF_{t,s_1}\right)\leq C_{\kappa}\|f\|_{L^q([s_1,s_2];L^p(\mR^d))},\label{Lp666}
\end{align}
where $C_{\kappa}$ is increasing with respect to $\kappa$, and $\sF_{t,s_1}$ is defined by (\ref{Flow}).
\et

The following lemma is taken from \cite[p. 1, Lemma 1.1]{Po}.
\bl\label{Le1}
Fix $t<0$. Let $\{\beta(s)\}_{s\in[t,0]}$ be a nonnegative measurable ($\sF_{t,s}$)-adapted process.
Assume that for all $t\leq s_1\leq s_2\leq 0$,
$$
\mE\left(\int^{s_2}_{s_1}\beta(r)\dif r\Bigg|_{\sF_{t,s_1}}\right)\leq\rho(s_1,s_2),
$$
where $\rho(s_1,s_2)$ is a nonrandom interval function satisfying the following conditions:
\begin{enumerate}[(i)]
\item $\rho(s_1,s_2)\leq\rho(s_3,s_4)$ if $(s_1,s_2)\subset(s_3,s_4)$;
\item $\lim_{\delta\downarrow 0}\sup_{t\leq s_1\leq s_2\leq 0, |s_1-s_2|\leq \delta}\rho(s_1,s_2)=0$.
\end{enumerate}
Then for any $\gamma>0$,
$$
\mE\exp\left\{\gamma\int^0_t\beta(r)\dif r\right\}\leq 2^N,
$$
where $N\in\mN$ is chosen being such that for any $k=0,\cdots, N-1$,
$$
\rho(-(k+1)|t|/N,-k|t|/N)\leq \frac{1}{2\gamma}.
$$
\el
Let $f$ be a locally integrable function on $\mR^d$. The Hardy-Littlewood maximal function is defined by
$$
\cM f(x):=\sup_{r>0}\frac{1}{|B_r|}\int_{B_r} f(x+y)\dif y,
$$
where $B_r:=\{y\in\mR^d: |y|<r\}$ and $|B_r|$ is the volume of $B_r$.

We recall the following well known results (cf. \cite[Appdenix]{Re-Zh} and \cite[p. 5, Theorem 1]{St}).
\bl\label{Le0}
(i) For any $f\in \mW^{1,p}$, there exist $C_d>0$ and a null set $E$ such that for all $x,y\notin E$,
\begin{align}
|f(x)-f(y)|\leq {C_d}(\cM|\nabla{f}|(x)+\cM|\nabla{f}|(y))|x-y|.\label{Epp2}
\end{align}
(ii) For any $p>1$, there exists a constant $C_{d,p}>0$ such that for any $f\in L^p(\mR^d)$,
\begin{align}
\|\cM f\|_p\leq{C_{d,p}}\|f\|_p.\label{Epp1}
\end{align}
\el

Using the above three tools, we can derive the following important estimates for later use.
\bl\label{Le2}
Suppose  that {\bf (H)$_\alpha$} holds with $\alpha\in(1,2)$, and $p>\frac{2d}{\alpha}$. For any $U>0$, there exists a time
$T=T(U)\in[-1,0)$  independent of $\nu\geq 1$  such that for any
divergence free $u\in L^\infty([T,0]; \mW^{1,p}(\mR^d;\mR^d))$ with
\begin{align}
\sup_{t\in[T,0]}\|u_t\|_{1,p}\leq U,\label{PP2}
\end{align}
the unique solution $X_{t,s}(x)$ to  SDE (\ref{EL7}) belongs to $\cap_{\gamma\geq 1}\mW^{1,\gamma}_{loc}$ with respect to $x$,
and preserves the volume, and satisfies that for any $\gamma\geq 1$ and some $C=C(T,\gamma,U)>0$,
\begin{align}
\sup_{t\in[T,0]}\sup_{x\in\mR^d}\mE|\nabla X_{t,0}(x)|^\gamma\leq C,\label{Lp6}
\end{align}
and
\begin{align}
\sup_{t\in[T,0]}\sup_{x\in\mR^d}\mE|\nabla X_{t,0}(x)|^4\leq 2.\label{Lp7}
\end{align}
Moreover, for any $\varphi\in\mW^{1,p}(\mR^d;\mR^d)$, if we define
\begin{align}
w_t:=\bP\mE(\nabla^{\mathrm{t}} X_{t,0}\cdot (\varphi\circ X_{t,0})),\label{EL55}
\end{align}
then $w\in C([T,0];\mW^{1,p})$ and
\begin{align}
\p_i w_t=\bP\mE[\nabla^{\mathrm{t}} X_{t,0}\cdot(\nabla \varphi-\nabla^{\mathrm{t}} \varphi)\circ X_{t,0}\cdot\p_i X_{t,0}].\label{EL6}
\end{align}
\el
\begin{proof}
Under (\ref{PP2}), it has been proven in \cite[Theorem 1.1]{Zh3} (see also \cite{Pr}) that SDE (\ref{EL7}) admits a unique strong solution $X_{t,s}(x)$
for each $x\in\mR^d$. Since $u$ is divergence free, $x\mapsto X_{t,x}(x)$ preserves the volume.
Let $u^\eps_t(x):=u_t*\rho_\eps(x)$ be the mollifying approximation of $u$, where $(\rho_\eps)_{\eps\in(0,1)}$
is a family of mollifiers. Let $X^\eps_{t,s}(x)$ solve the following SDE
$$
X^\eps_{t,s}(x)=x+\int^s_t u^\eps_r(X^\eps_{t,r}(x))\dif r+\nu^{1/\alpha}(L_s-L_t), \ \ t\leq s\leq 0.
$$
Then
$$
\nabla X^\eps_{t,s}(x)=\mI+\int^s_t \nabla u^\eps_r(X^\eps_{t,r}(x))\cdot \nabla X^\eps_{t,r}(x)\dif r,
$$
and
$$
|\nabla X^\eps_{t,s}(x)|\leq 1
+\int^s_t |\nabla u^\eps_r(X^\eps_{t,r}(x))|\cdot |\nabla X^\eps_{t,r}(x)|\dif r,
$$
where $|\cdot|$ denotes the Hilbert-Schmidt norm for a matrix. By Gronwall's inequality,
\begin{align}
|\nabla X^\eps_{t,s}(x)|\leq \exp\left\{\int^s_{t}|\nabla u^\eps_r(X^\eps_{t,r}(x))|\dif r\right\}.\label{Lp8}
\end{align}
By Theorem \ref{Kry1}, one has that for any $q>\frac{p\alpha}{p\alpha-d}$ and all $-1\leq t\leq s_1\leq s_2\leq 0$,
\begin{align*}
\mE\left(\int^{s_2}_{s_1}|\nabla u^\eps_r(X^\eps_{t,r}(x))|\dif r\Big|_{\sF_{t,s_1}}\right)
&\leq C_{\|u^\eps\|_\infty} \|\nabla u^\eps\|_{L^q([s_1,s_2]; L^p)}\\
&\leq C_{\|u\|_{L^\infty([t,0];\mW^{1,p})}} \|u\|_{L^q([s_1,s_2]; \mW^{1,p})}\\
&\leq C_{U}U^q|s_2-s_1|^{1/q},
\end{align*}
where the second inequality is due to the Sobolev embedding relation $\mW^{1,p}\subset L^\infty$ provided that $p>d$.
Hence, by Lemma \ref{Le1}, for any $\gamma\geq 1$,
\begin{align}
\sup_{\eps\in(0,1)}\sup_{x\in\mR^d}\mE\exp\left\{\gamma\int^{0}_{t}
|\nabla u^\eps_r(X^\eps_{t,r}(x))|\dif r\right\}<+\infty,\label{Lp1}
\end{align}
and one can choose $T\in[-1,0)$ depending on $U$ such that for all $t\in[T,0)$,
\begin{align}
\sup_{\eps\in(0,1)}\sup_{x\in\mR^d}\mE\exp\left\{4\int^{0}_{t}
|\nabla u^\eps_r(X^\eps_{t,r}(x))|\dif r\right\}\leq 2.\label{Lp99}
\end{align}

Let $A_{t,s}(x)$ solve the following linear random ODE:
$$
A_{t,s}(x)=\mI+\int^s_t\nabla u_r(X_{t,r}(x))\cdot A_{t,r}(x)\dif r.
$$
{\it(Claim)}: For any $\delta\in[1,2)$,
\begin{align}
\lim_{\eps\to 0}\sup_{t\leq s;t,s\in[T,0]}\sup_{x\in\mR^d}\mE|X^\eps_{t,s}(x)-X_{t,s}(x)|^\delta=0,\label{Lp3}
\end{align}
and for each $t\in[T,0)$ and $x\in\mR^d$,
\begin{align}
\lim_{\eps\to 0}\mE|\nabla X^\eps_{t,0}(x)-A_{t,0}(x)|^\delta=0.\label{Lp4}
\end{align}

We first prove (\ref{Lp3}). By (\ref{Epp2}), we have
\begin{align*}
|X^\eps_{t,s}(x)-X_{t,s}(x)|&\leq \int^s_t|u^\eps_r(X^\eps_{t,r}(x))-u^\eps_r(X_{t,r}(x))|\dif r
+\int^s_t|u^\eps_r(X_{t,r}(x))-u_r(X_{t,r}(x))|\dif r\\
&\leq C_d\int^s_t(\cM|\nabla u^\eps_r|(X^\eps_{t,r}(x))+\cM|\nabla u^\eps_r|(X_{t,r}(x)))|X^\eps_{t,r}(x)-X_{t,r}(x)|\dif r\\
&\quad+\int^0_t|u^\eps_r(X_{t,r}(x))-u_r(X_{t,r}(x))|\dif r,
\end{align*}
which yields by Gronwall's inequality that
\begin{align*}
|X^\eps_{t,s}(x)-X_{t,s}(x)|&\leq \exp\left\{C_d\int^s_t(\cM|\nabla u^\eps_r|(X^\eps_{t,r}(x))
+\cM|\nabla u^\eps_r|(\hat X_{t,r}(x)))\dif r\right\}\\
&\quad\times\int^0_t|u^\eps_r(X_{t,r}(x))-u_r(X_{t,r}(x))|\dif r.
\end{align*}
As in estimating (\ref{Lp1}), one has that for any $\gamma\geq 1$,
\begin{align*}
\sup_{\eps\in(0,1)}\sup_{x\in\mR^d}\mE\exp\left\{\gamma\int^{0}_{t}\cM|\nabla u^\eps_r|(X^\eps_{t,r}(x))\dif r\right\}<+\infty
\end{align*}
and
\begin{align*}
\sup_{\eps\in(0,1)}\sup_{x\in\mR^d}\mE\exp\left\{\gamma\int^{0}_{t}
\cM|\nabla u^\eps_r|(X_{t,r}(x))\dif r\right\}<+\infty.
\end{align*}
Hence, by H\"older's inequality and Theorem \ref{Kry1} again, we have for any $\delta\in[1,2)$ and $q>\frac{p\alpha}{p\alpha-2d}$,
\begin{align}
\mE|X^\eps_{t,s}(x)-X_{t,s}(x)|^\delta&\leq
\left(\mE\exp\left\{\frac{4\delta C_d}{2-\delta}\int^0_t\cM|\nabla u^\eps_r|(X^\eps_{t,r}(x))\dif r\right\}\right)^{(2-\delta)/4}\no\\
&\quad\times \left(\mE\exp\left\{\frac{4\delta C_d}{2-\delta}\int^0_t\cM|\nabla u^\eps_r|(X_{t,r}(x)))\dif r\right\}\right)^{(2-\delta)/4}\no\\
&\quad\times \left(|t|\mE\int^0_t|u^\eps_r-u_r|^2(X_{t,r}(x))\dif r\right)^{\delta/2}\no\\
&\leq C\||u^\eps-u|^2\|^{\delta/2}_{L^q([t,0];L^{p/2})}\no\\
&=C\|u^\eps-u\|^\delta_{L^{2q}([t,0];L^{p})},\label{PP3}
\end{align}
where $C$ is independent of $\eps,x,t,s$. Limit (\ref{Lp3}) now follows from the property of convolutions.

As for (\ref{Lp4}), we have
\begin{align*}
|\nabla X^\eps_{t,s}(x)-A_{t,s}(x)|&\leq\int^s_t |\nabla u^\eps_r(X^\eps_{t,r}(x))-\nabla u_r(X_{t,r}(x))|\cdot |\nabla X^\eps_{t,r}(x)|\dif r\\
&\quad+\int^s_t |\nabla u_r(X_{t,r}(x))|\cdot |\nabla X^\eps_{t,r}(x)-A_{t,r}(x)|\dif r,
\end{align*}
which then gives that
\begin{align*}
|\nabla X^\eps_{t,0}(x)-A_{t,0}(x)|&\leq\int^0_t |\nabla u^\eps_r(X^\eps_{t,r}(x))-\nabla u_r(X_{t,r}(x))|\cdot |\nabla X^\eps_{t,r}(x)|\dif r\\
&\quad\times\exp\left\{\int^{0}_{t}|\nabla u_r(X_{t,r}(x))|\dif r\right\}.
\end{align*}
As above, using (\ref{Lp8}), (\ref{Lp1}) and H\"older's inequality, we have
\begin{align}
\mE|\nabla X^\eps_{t,0}(x)-A_{t,0}(x)|^\delta&\leq C\left(\mE\int^0_t
|\nabla u^\eps_r(X^\eps_{t,r}(x))-\nabla u_r(X_{t,r}(x))|^2\dif r\right)^{\delta/2}.\label{Op1}
\end{align}
For fixed $\eps'\in(0,1)$, by (\ref{Lp3}), we have
\begin{align}
\mE\int^0_t|\nabla u^{\eps'}_r(X^\eps_{t,r}(x))-\nabla u^{\eps'}_r(X_{t,r}(x))|^2\dif r\to 0\ \mbox{ as } \eps\to 0.\label{Op2}
\end{align}
By (\ref{Lp666}), we have for $q>\frac{p\alpha}{p\alpha-2d}$,
\begin{align}
\sup_{\eps\in(0,1)}\mE\int^0_t|\nabla (u^{\eps'}_r-u_r)(X^\eps_{t,r}(x))|^2\dif r\leq
C\|\nabla (u^{\eps'}-u)\|_{L^{2q}([t,0];L^p)}^2\to 0\ \mbox{ as } \eps'\to 0.\label{Op3}
\end{align}
Limit (\ref{Lp4}) then follows by (\ref{Op1}), (\ref{Op2}) and (\ref{Op3}).

\vspace{5mm}

Using the above claims,  by (\ref{Lp8}), (\ref{Lp1}) and (\ref{Lp99}),
one finds that $X_{t,s}(\cdot)\in\cap_{\gamma\geq 1}\mW^{1,\gamma}_{loc}$, and (\ref{Lp6}) and (\ref{Lp7}) hold.
Moreover, if we define
$$
w^\eps_t:=\bP\mE(\nabla^{\mathrm{t}} X^\eps_{t,0}\cdot (\varphi\circ X^\eps_{t,0})),
$$
then, since $\bP$ is a bounded linear operator in $L^p$ and $x\mapsto X^\eps_{t,0}(x)$ preserves the volume,
it is easy to see that $w^\eps\in C([T,0];\cap_{k\in\mN}\mW^{k,p})$. Noting that
\begin{align}
0=\bP\nabla\mE( X^\eps_{t,0}\cdot (\varphi\circ X^\eps_{t,0}))=
\bP\mE(\nabla^{\mathrm{t}} X^\eps_{t,0}\cdot (\varphi\circ X^\eps_{t,0}))
+\bP\mE(\nabla^{\mathrm{t}} (\varphi\circ X^\eps_{t,0})\cdot X^\eps_{t,0}),\label{PP4}
\end{align}
we have
$$
\p_i w^\eps_t=\bP\mE[\nabla^{\mathrm{t}} X^\eps_{t,0}\cdot (\nabla \varphi-\nabla^{\mathrm{t}} \varphi)\circ X^\eps_{t,0}\cdot\p_i X^\eps_{t,0}].
$$
Using limits (\ref{Lp3}) and (\ref{Lp4}), formula (\ref{EL6}) then follows, and meanwhile,
$$
\lim_{\eps\to 0}\sup_{t\in[T,0]}\|w^\eps_t-w_t\|_{1,p}=0.
$$
So, $w\in C([T,0];\mW^{1,p})$. The proof is complete.
\end{proof}

The following lemma gives the continuous dependence of solutions to SDE (\ref{EL7}) with respect to $u$.
\bl\label{Le3}
Suppose  that {\bf (H)$_\alpha$} holds with $\alpha\in(1,2)$, and $p>\frac{2d}{\alpha}$. For  $U>0$
and $T\in[-1,0)$, let $u, \hat u\in L^\infty([T,0];\mW^{1,p}(\mR^d;\mR^d))$ be divergence free with
$$
\sup_{t\in[T,0]}\|u_t\|_{1,p}\leq U,\ \  \sup_{t\in[T,0]}\|\hat u_t\|_{1,p}\leq U.
$$
Let $X,\hat X$ be the solutions of SDE (\ref{EL7}) corresponding to $u,\hat u$.
Then for any $\delta\in[1,2)$, $q>\frac{p\alpha}{p\alpha-2d}$ and $t\in[T,0]$,
\begin{align}
\sup_{x\in\mR^d}\mE|X_{t,0}(x)-\hat X_{t,0}(x)|^\delta\leq C_1\|u-\hat u\|^\delta_{L^{2q}([t,0]; L^p)},\label{Op4}
\end{align}
where $C_1$ only depends on $U, T,\alpha,\delta,p,q,d,\psi$.
Moreover, for any $\varphi,\hat \varphi\in\mW^{1,p}(\mR^d;\mR^d)$, let $w_t$ and $\hat w_t$ be defined as in (\ref{EL55}) corresponding to $(\varphi,X)$ and
$(\hat \varphi,\hat X)$, then
for any $q>\frac{p\alpha}{p\alpha-2d}$ and $t\in[T,0)$,
\begin{align}
\|w_t-\hat w_t\|^{2q}_p\leq C_2\|\varphi-\hat \varphi\|^{2q}_p+C_3\int^0_t\|u_r-\hat u_r\|^{2q}_p\dif r,\label{Op5}
\end{align}
where $C_2$ (resp. $C_3$) only depends on $U, T,\alpha,p,q,d,\psi$ (resp. $U, T,\alpha,p,q,d,\psi,\|\varphi\|_{1,p}$).
\el
\begin{proof}
Estimate (\ref{Op4}) follows from the same calculations as in estimating (\ref{PP3}). Let us look at (\ref{Op5}).
Using mollifying approximation and (\ref{PP4}), we have
\begin{align*}
w_t-\hat w_t&=\bP\mE[\nabla^{\mathrm{t}} (X_{t,0}-\hat X_{t,0})\cdot (\varphi\circ X_{t,0})]\\
&\quad+\bP\mE[\nabla^{\mathrm{t}} \hat X_{t,0}\cdot (\varphi\circ X_{t,0}-\hat \varphi\circ\hat X_{t,0})]\\
&=\bP\mE[ \nabla^{\mathrm{t}}(\varphi\circ X_{t,0})\cdot(X_{t,0}-\hat X_{t,0})]\\
&\quad+\bP\mE[\nabla^{\mathrm{t}} \hat X_{t,0}\cdot (\varphi\circ X_{t,0}-\varphi\circ\hat X_{t,0})]\\
&\quad+\bP\mE[\nabla^{\mathrm{t}} \hat X_{t,0}\cdot (\varphi\circ \hat X_{t,0}-\hat \varphi\circ\hat X_{t,0})].
\end{align*}
By the boundedness of $\bP$ in $L^p$, we have
\begin{align*}
\|w_t-\hat w_t\|_p&\leq C\|\mE[\nabla^{\mathrm{t}}(\varphi\circ X_{t,0})\cdot (X_{t,0}-\hat X_{t,0})]\|_p\\
&\quad+C\|\mE[\nabla^{\mathrm{t}} \hat X_{t,0}\cdot (\varphi\circ X_{t,0}-\varphi\circ \hat X_{t,0})]\|_p\\
&\quad+C\|\mE[\nabla^{\mathrm{t}} \hat X_{t,0}\cdot (\varphi\circ\hat  X_{t,0}-\hat \varphi\circ \hat X_{t,0})]\|_p.
\end{align*}
By H\"older's inequality, for any $\delta\in(p/(p-1),2)$ and some $\beta=\beta(\delta,p)>1$, we have
\begin{align*}
&\|\mE[\nabla^{\mathrm{t}} \hat X_{t,0}\cdot (\varphi\circ X_{t,0}-\varphi\circ \hat X_{t,0})]\|_p\\
&\quad\stackrel{(\ref{Epp2})}{\leq} C\|\mE[|\nabla\hat X_{t,0}|\cdot (\cM|\nabla \varphi|(X_{t,0})+\cM|\nabla \varphi|(\hat X_{t,0}))\cdot
|X_{t,0}-\hat X_{t,0}|]\|_p\\
&\quad\leq C\|\|\nabla \hat X_{t,0}\|_{L^{\beta}(\Omega)}\cdot \|\cM|\nabla \varphi|(X_{t,0})+\cM|\nabla \varphi|(\hat X_{t,0})\|_{L^p(\Omega)}\cdot
\|X_{t,0}-\hat X_{t,0}\|_{L^\delta(\Omega)}\|_p\\
&\quad\leq C\sup_{x\in\mR^d}\|\nabla \hat X_{t,0}(x)\|_{L^{\beta}(\Omega)}\cdot
\sup_{x\in\mR^d}\|X_{t,0}(x)-\hat X_{t,0}(x)\|_{L^\delta(\Omega)}\cdot \|\cM|\nabla \varphi|\|_p\\
&\quad\stackrel{(\ref{Lp6})(\ref{Op4})(\ref{Epp1})}{\leq}  C\|u-\hat u\|_{L^{2q}([t,0]; L^p)}\cdot\|\nabla \varphi\|_p,
\end{align*}
where we have used that $x\mapsto X_{t,0}(x),\hat X_{t,0}(x)$ preserve the volume.
Similarly, we have
$$
\|\mE[\nabla^{\mathrm{t}}(\varphi\circ X_{t,0})\cdot (X_{t,0}-\hat X_{t,0})]\|_p
\leq C\|u-\hat u\|_{L^{2q}([t,0]; L^p)}\cdot\|\nabla \varphi\|_p,
$$
and
$$
\|\mE[\nabla^{\mathrm{t}} \hat X_{t,0}\cdot (\varphi\circ \hat X_{t,0}-\hat \varphi\circ \hat X_{t,0})]\|_p
\leq C\|\varphi-\hat \varphi\|_p.
$$
The proof is thus complete.
\end{proof}

We are now in a position to prove the following main result of this section.
\bt\label{Main1}
For any divergence free $u_0\in \mW^{1,p}(\mR^d;\mR^d)$ with $p>\frac{2d}{\alpha}$, there exist a time $T=T(\|u_0\|_{1,p})<0$ independent of $\nu\geq 1$ and
a unique pair of $(u,X)$ with $u\in C([T,0],\mW^{1,p})$ solving the following stochastic system:
\begin{align}\label{PP8}
\left\{
\begin{aligned}
X_{t,s}(x)&=x+\int^s_t u_r(X_{t,r}(x))\dif r+\nu^{1/\alpha}(L_s-L_t), t\leq s\leq 0,\\
u_t&=\bP\mE[\nabla^{\mathrm{t}} X_{t,0}\cdot (u_0\circ X_{t,0})].
\end{aligned}
\right.
\end{align}
Moreover, for some $C_0\geq 1$ only depending on $p$,
$$
\sup_{t\in[T,0]}\|u_t\|_{1,p}\leq 3C_0\|u_0\|_{1,p},
$$
and $u$ satisfies (\ref{FNS}) in a generalized sense, i.e., for all divergence free vector field $\phi\in C^\infty_0(\mR^d;\mR^d)$,
\begin{align}
\<u_t,\phi\>=\<u_0,\phi\>+\int^0_t[\<u_s,\cL^*_\nu\phi\>+\<(u_s\cdot\nabla)u_s,\phi\>]\dif s,\label{Eq4}
\end{align}
where
$$
\cL^*_\nu g(x)=\int_{\mR^d} \Big[g(x-\nu^{1/\alpha}y)-g(x)+(\nu^{1/\alpha}y\cdot \nabla )g(x)1_{|y|\leq 1}\Big]\nu(\dif y).
$$
\et
\begin{proof}
Set $u^0_r(x)\equiv u_0(x)$. For $n\in\mN_0$, by Lemma \ref{Le2}, we  recursively define
\begin{align}
\label{PP7}\left\{
\begin{aligned}
X^{n}_{t,s}(x)&=x+\int^s_t u^n_r(X^n_{t,r}(x))\dif r+\nu^{1/\alpha}(L_s-L_t),\\
u^{n+1}_t&=\bP\mE[\nabla^{\mathrm{t}} X^n_{t,0}\cdot (u_0\circ X^n_{t,0})].
\end{aligned}
\right.
\end{align}
Let us estimate the $\mW^{1,p}$-norm of $u^{n+1}$. First of all, by H\"older's inequality, we have
\begin{align*}
\|u^{n+1}_t\|_p&\leq C\|\mE[\nabla^{\mathrm{t}} X^n_{t,0}\cdot (u_0\circ X^n_{t,0})]\|_p\\
&\leq C\left(\int_{\mR^d}\|\nabla X^n_{t,0}(x)\|^p_{L^2(\Omega)}\cdot \|u_0\circ X^n_{t,0}(x)\|^p_{L^2(\Omega)}\dif x\right)^{1/p}\\
&\leq C\sup_{x\in\mR^d}\|\nabla X^n_{t,0}(x)\|_{L^2(\Omega)}
\left(\int_{\mR^d}\mE|u_0\circ X^n_{t,0}(x)|^p\dif x\right)^{1/p}\\
&=C\sup_{x\in\mR^d}\|\nabla X^n_{t,0}(x)\|_{L^2(\Omega)}\|u_0\|_p,
\end{align*}
 and by (\ref{EL6}),
\begin{align*}
\|\nabla u^{n+1}_t\|_p\leq C\|\mE[|\nabla X^n_{t,0}|^2\cdot |\nabla u_0\circ X^n_{t,0}|]\|_p
\leq C\sup_{x\in\mR^d}\|\nabla X^n_{t,0}(x)\|^2_{L^4(\Omega)}\|\nabla u_0\|_p,
\end{align*}
Hence, for some $C_0\geq 1$ only depending on $p$,
$$
\|u^{n+1}_t\|_{1,p}\leq C_0\left(1+\sup_{x\in\mR^d}\|\nabla X^n_{t,0}(x)\|^4_{L^4(\Omega)}\right)\|u_0\|_{1,p}.
$$

Now, taking $U=3C_0\|u_0\|_{1,p}$ in Lemma \ref{Le2}, by induction method, there exists a
time $T=T(U)<0$ independent of $\nu\geq 1$ such that for all $n\in\mN_0$,
\begin{align}
\sup_{t\in[T,0]}\|u^{n+1}_t\|_{1,p}\leq U,\quad \sup_{t\in[T,0]}\sup_{x\in\mR^d}\mE|\nabla X^n_{t,0}(x)|^4\leq 2,\label{PP5}
\end{align}
and for any $\gamma\geq 1$,
$$
\sup_{n\in\mN_0}\sup_{t\in[T,0]}\sup_{x\in\mR^d}\mE|\nabla X^n_{t,0}(x)|^\gamma<+\infty.
$$
Thus, by Lemma \ref{Le3}, one has that for all $t\in[T,0]$,
$$
\|u^{n+1}_t-u^{m+1}_t\|^{2q}_p\leq C\int^0_t\|u^n_r-u^m_r\|^{2q}_p\dif r,
$$
where the constant $C$ is independent of $n,m$ and $t$. By Gronwall's inequality, we get
$$
\lim_{n,m\to\infty}\sup_{t\in[T,0]}\|u^n_t-u^m_t\|^{2q}_p=0.
$$
Thus, by (\ref{PP5}), there exists a $u\in L^\infty([T,0];\mW^{1,p})$  such that
$$
\lim_{n\to\infty}\sup_{t\in[T,0]}\|u^n_t-u_t\|_{p}=0
$$
and
$$
\sup_{t\in[T,0]}\|u_t\|_{1,p}\leq U.
$$
Let $X$ be the solution of SDE (\ref{EL7}) corresponding to the above $u$. By taking limits for both sides of (\ref{PP7})
and using Lemma \ref{Le3}, one finds that
$(u,X)$ solves (\ref{PP8}). Moreover, the regularity of $u\in C([T,0];\mW^{1,p})$ follows from Lemma \ref{Le2}, and
the uniqueness follows from Lemma \ref{Le3}. As for equation (\ref{Eq4}), let $u^\eps_0(x):=u_0*\rho_\eps(x)$ be the mollifying approximation of $u_0$.
Let $u^\eps_t(x)\in C([T,0];\cap_{k\geq 0}\mW^{k,p})$ be the corresponding solution of (\ref{PP8}). By Theorem \ref{Rep},
one has that for all divergence free vector field $\phi\in C^\infty_0(\mR^d;\mR^d)$
$$
\<u^\eps_t,\phi\>=\<u^\eps_0,\phi\>+\int^0_t[\<u^\eps_s,\cL^*_\nu\phi\>+\<(u^\eps_s\cdot\nabla)u^\eps_s,\phi\>]\dif s.
$$
Taking limits $\eps\to 0$ and by Lemma \ref{Le3}, one then obtains equation (\ref{Eq4}).
\end{proof}

\section{Gradient estimates for L\'evy processes with drifts}

In this section we prove gradient estimates for SDE (\ref{EL7}). It will be used to obtain the global well posedness for stochastic system
(\ref{PP8}) with large viscosity $\nu$ in the periodic case. This section can be read independently and has some interests in itself.

Let $\cP(\mR^d)$ be the space of all Borel probability measures on $\mR^d$. For $\mu\in\cP(\mR^d)$, define
$$
\cT_\mu f(x)=\int_{\mR^d}f(x+y)\mu(\dif y),\ f\in C_b(\mR^d).
$$
It is clear that for any $\mu,\nu\in\cP(\mR^d)$,
$$
\cT_\mu\cT_\nu f=\cT_{\mu*\nu} f,\ f\in C_b(\mR^d),
$$
where $\mu*\nu$ denotes the convolution between $\mu$ and $\nu$.

For $a>0$, let us consider the truncated symbol of $\psi(\xi)$ in (\ref{Lp9}):
$$
\psi_a(\xi)=\int_{0<|x|\leq a}(1-e^{i\xi\cdot x}+i\xi\cdot x 1_{|x|\leq 1})\nu(\dif x).
$$
For any $t<0$, by Bochner's theorem, there exists a unique probability measure  $\mu_{a,t}\in\cP(\mR^d)$ such that
\begin{align}
\int_{\mR^d}e^{i\xi\cdot x}\mu_{a,t}(\dif x)=e^{t\psi_a(\xi)},\label{EL1}
\end{align}
and respectively, $\tilde\mu_{a,t}\in\cP(\mR^d)$ such that
\begin{align}
\int_{\mR^d}e^{i\xi\cdot x}\tilde\mu_{a,t}(\dif x)=e^{t(\psi(\xi)-\psi_a(\xi))}.\label{EL2}
\end{align}

We recall the following result from \cite[Proposition 2.3]{Sc-Sz-Wa}.
\bp\label{Pr1}
Assume that {\bf (H)$_\alpha$} holds with $\alpha\in(0,2)$. Then for any $a,t>0$,
$$
\mu_{a,t}(\dif x)=p_{a,t}(x)\dif x,
$$
and for any $n\in\mN_0$, there exists a time $t_0=t_0(n,d,\alpha,\psi)<0$ and a constant
$C=C(t_0,n,d,\alpha,\psi)>0$ such that for any $t\in[t_0,0)$,
\begin{align}
|\nabla^n p_{|t|^{1/\alpha},t}(x)|\leq C |t|^{-(d+n)/\alpha}(1+|t|^{-1/\alpha}|x|)^{-d-1}.\label{PL2}
\end{align}
\ep

For $\beta>0$, let $\mW^{\beta,p}:=(I-\Delta)^{-\beta/2}(L^p)$ be the Bessel potential space with the norm:
$$
\|f\|_{\beta,p}:=\|(I-\Delta)^{\beta/2}f\|_p.
$$
If $\beta=k\in\mN$, $\mW^{\beta,p}$ is the same as the Sobolev space $\mW^{k,p}$.
Using Proposition \ref{Pr1}, we now derive the following useful result.
\bt\label{Th1}
Assume that {\bf (H)$_\alpha$} holds with $\alpha\in(0,2)$.
For $\nu>0$, define
\begin{align}
\cT^\nu_tf(x):=\mE(f(\nu^{\frac{1}{\alpha}} L_t+x)).\label{EL5}
\end{align}
Then for any $p\in[1,\infty]$ and $m,n\in\mN$, there exists a constant $C=C(p,m,n,d,\psi)$ such that for all $f\in\mW^{m,p}$,
\begin{align}
\|\cT^\nu_t f\|_{m+n,p}\leq C[\nu(|t|\wedge 1)]^{-m/\alpha}\|f\|_{n,p},\ \
\forall \nu>0, \forall t<0.\label{Pl1}
\end{align}
If $p\in(1,\infty)$, the above estimate also holds for any nonnegative real numbers $m,n$.
\et
\begin{proof}
For $\nu,a>0$, let $\mu^\nu_{a,t}, \tilde \mu^\nu_{a,t}\in\cP(\mR^d)$
be defined as in (\ref{EL1}) and (\ref{EL2}) corresponding to the symbols $\psi_a(\nu^{1/\alpha}\xi)$ and
$\psi(\nu^{1/\alpha}\xi)-\psi_a(\nu^{1/\alpha}\xi)$.
It is easy to see that
\begin{align}
\mu^\nu_{a,t}(\dif x)=\nu^{-\frac{d}{\alpha}} p_{a,t}(\nu^{-\frac{1}{\alpha}}x)\dif x.\label{EL3}
\end{align}
Set
$$
\cT^\nu_{a,t}f=\cT_{\mu^\nu_{a,t}}f,\ \ \tilde\cT^\nu_{a,t}f=\cT_{\tilde\mu^\nu_{a,t}}f,
$$
then for all $t\in\mR_-$,
\begin{align}
\cT^\nu_tf=\cT^\nu_{a,t}\tilde\cT^\nu_{a,t}f.\label{EL4}
\end{align}
By (\ref{EL3}) and Proposition \ref{Pr1}, there exists a time $t_0=t_0(n,d,\alpha,\psi)<0$ such that for all $t\in[t_0,0)$ and $f\in L^p(\mR^d)$,
\begin{align*}
\nabla^n\cT^\nu_{|t|^{1/\alpha},t}f(x)&=\nu^{-\frac{d+n}{\alpha}}
\int_{\mR^d}f(y)\nabla^np_{|t|^{1/\alpha},t}(\nu^{-\frac{1}{\alpha}}(y-x))\dif y\\
&=\nu^{-\frac{d+n}{\alpha}}
\int_{\mR^d}f(x+y)\nabla^np_{|t|^{1/\alpha},t}(\nu^{-\frac{1}{\alpha}}y)\dif y.
\end{align*}
Hence, for any $p\in[1,\infty]$,  by Minkowskii's inequality and (\ref{PL2}), we have
\begin{align}
\|\nabla^n\cT^\nu_{|t|^{1/\alpha},t}f\|_p&\leq\nu^{-\frac{d+n}{\alpha}}
\|f\|_p\int_{\mR^d} |\nabla^np_{|t|^{1/\alpha},t}(\nu^{-\frac{1}{\alpha}}y)|\dif y\no\\
&\leq C(\nu |t|)^{-(d+n)/\alpha}\|f\|_p\int_{\mR^d} (1+(\nu |t|)^{-1/\alpha}|y|)^{-d-1}\dif y\no\\
&= C(\nu |t|)^{-n/\alpha}\|f\|_p\int_{\mR^d} (1+|y|)^{-d-1}\dif y=\tilde C(\nu |t|)^{-n/\alpha}\|f\|_p.
\end{align}
Thus, by (\ref{EL4}) and the $L^p$-contraction of $\tilde \cT_{a,t}$, we have
\begin{align*}
\|\nabla^{n+m}\cT^\nu_tf\|_p\leq \tilde C(\nu|t|)^{-n/\alpha}\|\tilde\cT^\nu_{|t|^{1/\alpha},t}\nabla^mf\|_p\leq
\tilde C(\nu|t|)^{-n/\alpha}\|\nabla^mf\|_p,
\end{align*}
which yields that
$$
\|\cT^\nu_tf\|_{m+n,p}\leq \tilde C(\nu|t|)^{-n/\alpha}\|f\|_{m,p}.
$$
For general $t<t_0$, it follows by the semigroup property of $\cT^\nu_t$.
If $p\in(1,\infty)$ and $m,n$ are nonnegative real numbers, it follows by interpolation theorem (cf. \cite{Tr}).
\end{proof}

Let us consider the following PIDE:
\begin{align}
\p_t h+\cL_\nu h+(u\cdot\nabla)h=0,\ \ t\leq 0,\label{PP6}
\end{align}
subject to the final value
$$
h_0(x)=\varphi(x),
$$
where
$$
\cL_\nu u(x):=\int_{\mR^d} \Big[u(x+\nu^{1/\alpha}y)-u(x)-(\nu^{1/\alpha}y\cdot \nabla )u(x)1_{|y|\leq 1}\Big]\nu(\dif y).
$$
By Duhamel's formula, one can write equation (\ref{PP6}) as the following mild form:
\begin{align}
h_t(x)=\cT^\nu_t\varphi(x)+\int^0_t\cT^\nu_{t-s}((u_s\cdot\nabla) h_s)(x)\dif s,\label{PP9}
\end{align}
where $\cT^\nu_t$ is defined by (\ref{EL5}).

We need the following simple lemma (cf. \cite[Lemma 5.1]{Zh1}).
\bl\label{l2}
Let $z(t)$  be a nonnegative function defined on $[-1,0)$, and satisfy that for some $K_1, K_2>0$ and $\b,\gamma\in(0,1)$,
$$
z(t)\leq K_1 |t|^{-\gamma}+K_2\int^0_t\frac{z(s)}{(s-t)^\beta}\dif s,\quad\forall t\in[-1,0).
$$
Then for any $t\in[-1,0)$
$$
z(t)\leq C_{K_2,\b,\gamma}K_1|t|^{-\gamma}.
$$
\el

We have the following useful estimate.
\bl\label{Le4}
Assume that {\bf (H)$_\alpha$} holds with $\alpha\in(1,2)$ and
$$
u\in L^\infty([-1,0]; L^p(\mR^d)+L^\infty(\mR^d))
$$
provided that $p>\frac{d}{\alpha-1}$. Then for any $\gamma\in(\frac{d}{p}+1,\alpha)$ and $\varphi\in L^p(\mR^d)$,
there exists a unique solution $h\in L^1([-1,0]; \mW^{\gamma,p})$ to equation  (\ref{PP9})
such that for all $t\in[-1,0)$ and $\nu\geq 1$,
\begin{align}
\|h_t\|_{\gamma,p}\leq C_1(\nu|t|)^{-\gamma/\alpha}\|\varphi\|_p,\label{Lp2}
\end{align}
where $C_1$ only depends on $\gamma,p,d,\alpha,\psi$ and the norm of $\|u\|_{L^\infty([-1,0]; L^p(\mR^d)+L^\infty(\mR^d))}$.
In the case of $u\in L^\infty([-1,0]\times\mR^d)$,  for any $p\in[1,\infty]$ and $\varphi\in L^p(\mR^d)$,
there exists a unique solution $h\in L^1([-1,0]; \mW^{1,p})$ to equation (\ref{PP9}) such that  for all $t\in[-1,0)$  and $\nu\geq 1$,
\begin{align}
\|\nabla h_t\|_p\leq C_2(\nu|t|)^{-1/\alpha}\|\varphi\|_p,\label{PP1}
\end{align}
where $C_2$ only depends on $p,d,\alpha,\psi$ and the bound of $u$,  and is increasing with respect to the bound of $u$.
\el
\begin{proof}
We only prove the a priori estimates (\ref{Lp2}) and (\ref{PP1}). As for the existence, it follows from a standard Picardi iteration argument.
Assume that $u=u_1+u_2$ with $u_1\in L^\infty([-1,0]; L^p(\mR^d))$ and $u_2\in L^\infty([-1,0]; L^\infty(\mR^d))$. Then,
by (\ref{Pl1}) and the Sobolev embedding $\mW^{\gamma,p}\subset C^1_b$ provided $(\gamma-1) p>d$, we have
\begin{align*}
\|h_t\|_{\gamma,p}&\leq \|\cT^\nu_t\varphi\|_{\gamma,p}+C\int^0_t(\nu(s-t))^{-\frac{\gamma}{\alpha}}\|(u_s\cdot\nabla) h_s\|_p\dif s\\
&\leq \|\cT^\nu_t\varphi\|_{\gamma,p}+C\int^0_t(s-t)^{-\frac{\gamma}{\alpha}}
\Big(\|u_{1,s}\|_\infty\cdot\|\nabla h_s\|_p+\|u_{2,s}\|_p\cdot\|\nabla h_s\|_\infty\Big)\dif s\\
&\leq C(\nu|t|)^{-\gamma/\alpha}\|\varphi\|_{p}+C\int^0_t(s-t)^{-\frac{\gamma}{\alpha}}
\Big((\|u_{1,s}\|_\infty+\|u_{2,s}\|_p)\cdot\|\nabla h_s\|_{\gamma,p}\Big)\dif s.
\end{align*}
By Lemma \ref{l2}, we  obtain (\ref{Lp2}).

In the case of $u\in L^1([-,0];L^\infty(\mR^d))$, one has
\begin{align*}
\|\nabla h_t\|_p&\leq \|\nabla\cT^\nu_t\varphi\|_p+C\int^0_t(\nu(s-t))^{-\frac{1}{\alpha}}
\|(u_s\cdot\nabla)h_s\|_p\dif s\\
&\leq C(\nu|t|)^{-1/\alpha}\|\varphi\|_p+C\int^0_t(s-t)^{-\frac{1}{\alpha}}
\|u_s\|_\infty\cdot\|\nabla h_s\|_p\dif s,
\end{align*}
which gives (\ref{PP1}) by Lemma \ref{l2} again.
\end{proof}

Before proving the gradient estimates for SDE (\ref{EL7}), we recall the notions of weak existence and uniqueness for SDE (\ref{EL7}).
Fix $t<0$. Let $\mD_t$ be the space of all c\`adl\`ag functions from $[t,0]\to\mR^d$. We endow $\mD_t$ with the Skorohod metric so that
$\mD_t$ is a Polish space. A weak solution of SDE (\ref{EL7}) means that there exists a probability space $(\Omega,\sF,P)$ and two c\`adl\`ag processes $(X_{t,s})_{s\in[t,0]}$
and $(L_{t,s})_{s\in[t,0]}$ defined on it such that $s\mapsto L_{t,s}$ is a L\'evy process with symbol $\psi(\xi)$ and
$$
X_{t,s}(x)=x+\int^s_t u_r(X_{t,r}(x))\dif r+\nu^{1/\alpha} L_{t,s}, \ \ s\in[t,0],\ \ P-a.s.
$$
Such a solution will be denoted by $(\Omega,\sF,P; (X_{t,s})_{s\in[t,0]}, (L_{t,s})_{s\in[t,0]})$.
Weak uniqueness means that for two weak solutions $(\Omega^{(i)},\cF^{(i)},P^{(i)}; (X^{(i)}_{t,s})_{s\in[t,0]}, (L^{(i)}_{t,s})_{s\in[t,0]})$,
$i=1,2$, the laws of $s\mapsto X^{(i)}_{t,s}$ in $\mD_t$ are the same for $i=1,2$.

Assume that  for each $x\in\mR^d$, weak uniqueness holds for SDE (\ref{EL7}). Then for any bounded measurable
function $\varphi$, the mapping
$$
x\mapsto \mE \varphi(X_{t,0}(x))
$$
is well defined. Now, we can prove the following gradient estimate for $x\mapsto\mE \varphi(X_{t,0}(x))$.
\bt
Assume that {\bf (H)$_\alpha$} holds with $\alpha\in(1,2)$ and
\begin{align}
u\in L^\infty([-1,0]; L^p(\mR^d)+L^\infty(\mR^d))\label{EP1}
\end{align}
provided that $p>\frac{d}{\alpha-1}$. We also assume that for each $x\in\mR^d$, weak uniqueness holds for SDE (\ref{EL7}).
Then for any $\gamma\in(\frac{d}{p}+1,\alpha)$, $\nu\geq 1$ and $\varphi\in L^p(\mR^d)$,
\begin{align}
\|\mE\varphi(X_{t,0}(\cdot))\|_{\gamma,p}\leq C_1(\nu |t|)^{-\gamma/\alpha}\|\varphi\|_p;\label{Op66}
\end{align}
in the case of $u\in L^\infty([-1,0]\times\mR^d)$, for any $p\in[1,\infty]$, $\nu\geq 1$ and $\varphi\in L^p(\mR^d)$,
\begin{align}
\|\nabla\mE\varphi(X_{t,0}(\cdot))\|_p\leq C_2(\nu |t|)^{-1/\alpha}\|\varphi\|_p,\label{Op6}
\end{align}
where $C_1$ and $C_2$ are the same as in Lemma \ref{Le4}.
\et
\begin{proof}
Let $u^\eps_t(x)=u_t*\rho_\eps(x)$ be the mollifying approximation of $u_t$, and $X^\eps_{t,s}$
be the corresponding solution of SDE (\ref{EL7}). For $\varphi\in C^\infty_b(\mR^d)$, define
$$
h^\eps_t(x):=\mE\varphi(X^\eps_{t,0}(x)).
$$
By Theorem \ref{Th2}, $u^\eps$ solves the following PIDE:
$$
\p_t h^\eps+\cL_\nu h^\eps+(u^\eps\cdot\nabla)h^\eps=0,
$$
which is equivalent by Duhamel's formula that,
$$
h^\eps_t(x)=\cT^\nu_t\varphi(x)+\int^0_t\cT^\nu_{t-s}((u^\eps_s\cdot\nabla)h^\eps_s)(x)\dif s.
$$

Under (\ref{EP1}), it is well known that the laws of $\{X^\eps_{t,\cdot}(x)\}_{\eps\in(0,1)}$ in $\mD_t$ is tight and any accumulation point is a weak solution
of SDE (\ref{EL7}) (cf. \cite[Theorem 4.1]{Zh3}). By the weak uniqueness,
the whole sequence of $X^\eps_{t,\cdot}(x)$ weakly converges to the weak solution $X_{t,\cdot}(x)$ in $\mD_t$.
In particular, for any $x\in\mR^d$,
\begin{align}
\lim_{\eps\to 0}h^\eps_t(x)=\lim_{\eps\to 0}\mE\varphi(X^\eps_{t,0}(x))=\mE\varphi(X_{t,0}(x)).\label{EP3}
\end{align}
By (\ref{Lp2}), we have for all $\varphi\in C^\infty_0(\mR^d)$,
$$
\|\mE\varphi(X^\eps_{t,0}(\cdot))\|_{\gamma,p}=\|h^\eps_t\|_{\gamma,p}\leq C(\nu|t|)^{-\gamma/\alpha}\|\varphi\|_p.
$$
By (\ref{EP3}), we get (\ref{Op6}) for any $\varphi\in C^\infty_0(\mR^d)$. For general $\varphi\in L^p(\mR^d)$, it follows by a
standard approximation.

In the case of $u\in L^\infty([-1,0]\times\mR^d)$, by (\ref{PP1}), for $p=\infty$, we have
$$
|\mE\varphi(X^\eps_{t,0}(x))-\mE\varphi(X^\eps_{t,0}(y))|=|h^\eps_t(x)-h^\eps_t(y)|\leq C(\nu|t|)^{-1/\alpha}\|\varphi\|_\infty|x-y|.
$$
Letting $\eps\to 0$ yields that
$$
|\mE\varphi(X_{t,0}(x))-\mE\varphi(X_{t,0}(y))|\leq C(\nu|t|)^{-1/\alpha}\|\varphi\|_\infty|x-y|,
$$
which then gives  (\ref{Op6}) for $p=\infty$ and $\varphi\in C^\infty_b(\mR^d)$.
For general $\varphi\in L^\infty(\mR^d)$, it follows by a standard approximation.
For $p\in[1,\infty)$, it is similar.
\end{proof}

\section{Existence of global solutions for FNSE with large viscosity}

In this section, we prove the global existence of $\mW^{1,p}$-solution for stochastic system (\ref{PP8} ) in the case of  large viscosity and periodic boundary.
Let $\mT^d$ be the $d$-dimensional torus. Below, we shall work in the Sobolev space $\mW^{1,p}(\mT^d;\mR^d)$ with vanishing mean denoted by $\mW^{1,p}_0(\mT^d;\mR^d)$.
Thus, by Poinc\`are's inequality, for any $p>1$,
$$
\|u\|_p\leq C\|\nabla u\|_p,
$$
and an equivalent norm in $\mW^{1,p}(\mT^d;\mR^d)$ is thus given by
$$
\|u\|_{1,p}\simeq\|\nabla u\|_p.
$$
Below, we shall use $\|\nabla u\|_p$ as the norm of $\mW^{1,p}_0(\mT^d;\mR^d)$.
\bt
For any divergence free $u_0\in \mW^{1,p}_0(\mT^d;\mR^d)$ with $p>\frac{2d}{\alpha}$, there exist $\nu=\nu(\|u_0\|_{1,p})\geq 1$ sufficiently large
and a unique pair of $(u,X)$ with $u\in C(\mR_-,\mW^{1,p}_0(\mT^d;\mR^d))$ solving stochastic system (\ref{PP8}).
\et
\begin{proof}
By Theorem \ref{Main1}, there exist a time $T=T(\|\nabla u_0\|_{p})\in[-1,0)$ independent of $\nu\geq 1$ and a unique pair of $(u, X)$
with $u\in C([T,0],\mW^{1,p}_0)$ and
$$
\|\nabla u_t\|_p\leq 3C_0\|\nabla u_0\|_p
$$
solving  stochastic system (\ref{PP8}). Our aim is to prove that for some $T_*=T_*(\|\nabla u_0\|_{p})\in[T,0)$ and $\nu$ large enough,
\begin{align}
\|\nabla u_{T_*}\|_p\leq \|\nabla u_0\|_p.\label{EP4}
\end{align}
After proving this estimate,  one can invoke the standard bootstrap argument to obtain the global solution.

By (\ref{EL6}), we have
\begin{align*}
\p_iu_t&=\bP\mE[(\nabla^{\mathrm{t}} X_{t,0}-\mI)\cdot \nabla u_0\circ X_{t,0}\cdot\p_i X_{t,0}]+\bP\p_i\mE[u_0\circ X_{t,0}]\\
&+\bP\mE[\nabla^{\mathrm{t}} X_{t,0}\cdot \nabla^{\mathrm{t}} u_0\circ X_{t,0}\cdot\p_i(X_{t,0}-x)]+\bP\mE[\nabla(u_0^i\circ X_{t,0})].
\end{align*}
Noticing that
\begin{align*}
\nabla X_{t,s}(x)-\mI=\int^s_t\nabla u_r(X_{t,r}(x))\cdot\nabla X_{t,r}(x)\dif r,
\end{align*}
by H\"older's inequality, Theorem \ref{Kry1} and (\ref{Lp6}), we have for any $\gamma\in(1,\frac{p\alpha}{2d})$ and $x\in\mR^d$,
\begin{align*}
\mE|\nabla^{\mathrm{t}} X_{t,0}(x)-\mI|^2&\leq\mE\left[\int^0_t|\nabla u_r(X_{t,r}(x))|^2\dif r\int^0_t|\nabla X_{t,r}(x)|^2\dif r\right]\\
&\leq |t|^2\left[\mE\int^0_t|\nabla u_r(X_{t,r}(x))|^{2\gamma}\dif r\right]^{1/\gamma}\left[\mE\int^0_t|\nabla X_{t,r}(x)|^{2\gamma^*}
\dif r\right]^{1/\gamma^*}\leq C|t|^2,
\end{align*}
where $\gamma^*=\frac{\gamma}{\gamma-1}$ and $C$ is independent of $x$ and $t$, and depends on $\|\nabla u_0\|_p$.
Hence, by H\"older's inequality,
\begin{align*}
&\|\mE[(\nabla^{\mathrm{t}} X_{t,0}-\mI)\cdot \nabla u_0\circ X_{t,0}\cdot\p_i X_{t,0}]\|_p\\
&\quad\leq\|\|\nabla X_{t,0}-\mI\|_{L^2(\Omega)}\cdot \|\nabla u_0\circ X_{t,0}\|_{L^p(\Omega)}
\cdot\|\p_i X_{t,0}\|_{L^{2p/(p-2)}(\Omega)}\|_p\\
&\quad\leq C\sup_{x\in\mR^d}\|\nabla X_{t,0}(x)-\mI\|_{L^2(\Omega)}\cdot
\sup_{x\in\mR^d}\|\p_i X_{t,0}(x)\|_{L^{2p/(p-2)}(\Omega)}\cdot \|\nabla u_0\|_p\\
&\quad\leq C|t|\cdot \|\nabla u_0\|_p.
\end{align*}
Similarly,
$$
\|\mE[\nabla^{\mathrm{t}} X_{t,0}\cdot \nabla^{\mathrm{t}} u_0\circ X_{t,0}\cdot\p_i(X_{t,0}-x)]\|_p\leq C|t|\cdot \|\nabla u_0\|_p.
$$
On the other hand, by (\ref{Op6}), we have
$$
\|\p_i\mE[u_0\circ X_{t,0}]\|_p+\|\nabla\mE[u^i_0\circ X_{t,0}]\|_p\leq C(\nu|t|)^{-1/\alpha}\|u_0\|_p
\leq C(\nu|t|)^{-1/\alpha}\|\nabla u_0\|_p.
$$
Combining the above calculations, we obtain that
$$
\|\nabla u_t\|_p\leq C(|t|+(\nu|t|)^{-1/\alpha})\|\nabla u_0\|_p,
$$
which then gives (\ref{EP4}) by first letting $t$ small enough and then $\nu$ large enough. The proof is thus complete.
\end{proof}

{\bf Acknowledgements:}

The supports by NSFs of China (Nos. 10971076; 10871215) are  acknowledged.

\end{document}